\documentclass[12pt]{article}
\usepackage{amsmath,amsfonts,amsthm,amssymb}
\usepackage{bm}
\usepackage[dvips]{graphicx}

\newtheorem{Thm}{Theorem}[section]
\newtheorem{Def}{Definition}[section]
\newtheorem{Rem}{Remark}[section]
\newtheorem{Alg}{Algorithm}[section]

\newcommand{\df}{\stackrel{\mathrm{def.}}{=}}
\newcommand{\Rot}{\mathrm{Rot}}
\newcommand{\z}{{G\mbox{\textcircled{$\mathrm{z}$}}(H_1,H_2)}}
\newcommand{\zp}{{G\mbox{\textcircled{$\mathrm{z}$}}^\prime H}}
\newcommand{\Rotz}{\mathrm{Rot}_{G\mbox{\scriptsize{\textcircled{$\mathrm{z}$}}}(H_1,H_2)}}
\newcommand{\Rotzp}{\mathrm{Rot}_{G\mbox{\scriptsize{\textcircled{$\mathrm{z}$}}}^\prime H}}
\newcommand{\para}{{\scriptscriptstyle\parallel}}
\newcommand{\nullgraph}{\varnothing}

\title{Generalized zig-zag products of regular digraphs and bounds on their
  spectral expansions}
\author{Shunichi NOMURA and Akimichi TAKEMURA\\ \\
        Department of Mathematical Informatics\\
        Graduate School of Information Science and Technology\\
        The University of Tokyo}
\date{March, 2007}
\begin{document}

\maketitle

\begin{abstract}
  We introduce a generalization of the zig-zag product of regular
  digraphs (directed graphs), which allows us to construct regular digraphs with more
  flexible choices of the degrees.  In our
  generalization, we can control the connectivity of the resulting graph
  measured by its spectral expansion.
  We derive an upper bound on the spectral expansion of the
  generalized zig-zag product. Our upper bound improves on known
  bounds when applied to the zig-zag product. We also consider a
  special case of the generalized zig-zag product, where one of the
  components is a trivial graph whose edges are all self-loops.  We
  call it a reduced zig-zag product and derive a bound on the
  spectral expansion of its powers.
\end{abstract}

\section{Introduction}
A sparse graph which has high connectivity properties is called an
expander graph.  Expander graphs have many applications, 
such as complexity theory (\cite{Rei}), derandomization
(\cite{AKS}, \cite{AFWZ}) and error correcting codes (\cite{KSR},
\cite{SS}).  Connectivity properties of an expander graph are
measured by several expansion parameters, which are related to each
other.  In this paper we consider the spectral expansion as a measure
of connectivity. It can be interpreted as the rate at which a random walk
on the graph converges to its stationary distribution.

Many articles show that almost every regular graph has good expansion properties theoretically (\cite{Dur}, \cite{Fri})
or practically (\cite{HLW}, \cite{Novi}), based on randomization arguments.
Also many authors proposed explicit construction of expander graphs
using number theory or group theory
(\cite{Gra05}, \cite{LPS88}, \cite{Mar73}, \cite{Mar88}).
On the other hand, Reingold et al.\ \cite{RVW} introduced a new graph
operation called a zig-zag product and
succeeded in constructing  larger expander graphs by iteratively applying
the product.
The zig-zag product requires two graphs and generates a graph which
preserves the spectral expansions of the component graphs to some degree.
Originally the zig-zag product was defined only for regular undirected graphs,
but later Reingold et al.\ (\cite{RTV1}, \cite{RTV2}) extended the definition to directed graphs with an application to complexity theory.

In this paper, we introduce a generalization of the zig-zag product of
regular directed graphs.
It requires three graphs and also controls the resulting spectral
expansion in terms of the spectral expansions of the three components.
We derive an upper bound on the spectral expansion of the generalized 
zig-zag product. Our upper bound improves on the 
bounds in Reingold et al.\ (\cite{RTV1},\cite{RTV2}) when applied to the zig-zag
product.
Additionally we consider a special case of the 
generalized zig-zag product, where one of the components is 
a trivial graph whose edges are all self-loops.
We call it a reduced zig-zag product.
The spectral expansion of the reduced product itself is 1, which is
the worst case.  However we derive 
a bound for the powers of the reduced product, which is similar to the
bound in the generalized zig-zag product.

The organization of this paper is as follows.
In Section 2 we give notations and preliminary definitions.
In Section 3 we define the generalized zig-zag product and related graph operations.
In Section 4 we derive an upper bound for the spectral expansion 
of the zig-zag product and the powers of the reduced zig-zag product.
In Section 5 we present some results of numerical experiments to compare the spectral expansion of the zig-zag product
and its bound given in Section 4.

\section{Notation and Preliminaries}
In this paper, we consider digraphs 
which may have multiple edges and
self-loops.
In a digraph, the outdegree of a vertex is the number of edges leaving the vertex,
and the indegree of a vertex is the number of edges entering the vertex.
A digraph is $M$-outregular if every vertex has the outdegree $M$,
and $M$-inregular if every vertex has the indegree $M$.
A graph is $M$-regular if it is both $M$-outregular and $M$-inregular.

Given an $M$-regular graph $G$ on the set of vertices
$[N]\df\{1,\dots,N\}$,
consider a random walk on $G$ described by the transition matrix $A$
whose $(v,u)$'th entry is the number of edges from $u$ to $v$, divided by $M$\footnote
{Often the transition matrix is defined to be the transpose of our definition.
Following the recent literature on the zigzag product we adopt the 
the present definition of $A$. In our definition the transition of 
probability vector $\bm{x}$ by one step of the random work is
expressed by left-multiplication $A\bm{x}$ of $\bm{x}$ by $A$.}.
Let $\bm{1}_N = (1,\dots,1)^T \in \mathbb{R}^N$ denote the
$N$-dimensional vector consisting of $1$'s. 
By the regularity of $G$, $\bm{\pi} = \bm{1}_N/N$ is a stationary
distribution of the random walk, i.e. $A\bm{\pi} = \bm{\pi}$. 
Also note that $\bm{1}_N^T A =\bm{1}_N^T$ by the regularity of $A$.

We are interested in the rate at which random walks on $G$ converge to
a stationary distribution $\bm{\pi}$.
The convergence rate can be studied in terms of the 
Euclidean norm on $\mathbb{R}^N$.
Let $\langle \bm{x},\bm{y} \rangle = \sum_{v\in[N]}x_v y_v$ 
and $\Vert \bm{x}\Vert  = \sqrt{\langle \bm{x},\bm{x} \rangle}$
denote the standard inner product and the standard Euclidean norm
in $\mathbb{R}^N$.
We characterize the rate of convergence by the following parameter
called the spectral expansion.

\begin{Def}[spectral expansion]
Let $G$ be a regular digraph on the set of vertices $[N]$ and let $A$
be the transition matrix of the random walk on $G$.
The spectral expansion of $G$ is defined as
\[ \lambda(G) = \max_{\bm{x}\perp 1_N}\frac{\Vert A\bm{x}\Vert }{\Vert \bm{x}\Vert }. \]
\end{Def}

In the case that $G$ is undirected, $\lambda(G)$ is the second largest eigenvalue (in absolute value)
of the symmetric matrix $A$ (\cite{Bub}, \cite{Hag}).
In general, following Fill \cite{Fill} and Mihail \cite{Mih}, $\lambda(G)$ is
the square root of the second largest eigenvalue of $A^TA$,
which also means the second largest singular value of $A$.
If $\bm{\mu}$ is our initial distribution, 
then $\bm{\mu} - \bm{\pi}$ and $A(\bm{\mu} - \bm{\pi})$ are
orthogonal to $\bm{1}_N$.  Therefore
\[ \Vert A^t\bm{\mu} - \bm{\pi}\Vert  = \Vert A^t(\bm{\mu}-\bm{\pi})\Vert  \leq \lambda(G)^t \cdot \Vert \bm{\mu}-\bm{\pi}\Vert , \]
and the distance to $\bm{\pi}$ decreases exponentially at least with the rate 
of $\lambda(G)$.
In particular, if $\lambda(G) < 1$, then $\bm{\pi}$ is the unique stationary distribution.

In our analysis, we make use of 
the singular value decomposition of the transition matrix.
The singular value decomposition of $A$ is
\[ A = P\Sigma Q^T \]
where $P,Q$ are orthogonal matrices and $\Sigma$ is a diagonal matrix.
We denote $P=(\bm{p}_1,\dots,\bm{p}_N)$, 
$Q=(\bm{q}_1,\dots,\bm{q}_N)$ and $\Sigma=diag(\sigma_1,\dots,\sigma_N)$.
Then, for $v\in[N]$, $\bm{p}_v$ and $\bm{q}_v$ are the left-singular
and the right-singular vectors for $\sigma_v$,
respectively (i.e. $A\bm{q}_v=\sigma_v\bm{p}_v$ and $A^T\bm{p}_v=\sigma_v\bm{q}_v$).
By the regularity of $G$, $A\bm{1}_N=A^T\bm{1}_N=\bm{1}$, and we 
set $\sigma_1=1,\bm{p}_1=\bm{q}_1=\bm{1}_N/\sqrt{N}$.
Then, since $\bm{q}_2,\dots,\bm{q}_N\perp\bm{1}_N$, $\lambda(G)$ is 
the second largest singular value.

For vectors $\bm{x}\in\mathbb{R}^N$ and $\bm{y}\in\mathbb{R}^M$, we
define  their tensor product $\bm{x}\otimes \bm{y}$
to be the vector in $\mathbb{R}^{MN}$ whose $(u,k)$'th entry is $x_uy_k$.
Similarly, for an $N\times N$ matrix $A$ and an $M\times M$ matrix $B$, we define their tensor product $A\otimes B$
to be the $MN\times MN$ matrix whose $((v,l),(u,k))$'th entry is $A_{vu}B_{lk}$.
Then $(A\otimes B)(\bm{x}\otimes \bm{y}) = (A\bm{x})\otimes (B\bm{y})$.

\section{Operations on Directed Graphs}
In this section, we define some graph operations.
To define them, we give distinct labels to 
the edges leaving and entering each vertex by numbers from 1 to $M$,
which is called the two-way labelling.
Let $G$ be an $M$-regular digraph on the set of vertices $[N]$.
A two way labelling of $G$ is a family of bijections from 
$[M]$ to the set of edges leaving each vertex 
and the set of edges entering each vertex.
Such a graph together with its two-way labelling can be specified by a rotation map
$\Rot_G:[N]\times[M]\rightarrow[N]\times[M]$,
where $\Rot_G(u,k)=(v,l)$ if the $k$-th edge leaving $u$ is the $l$-th edge entering $v$.

First, we define a generalization of the zig-zag product introduced in \cite{RVW}.
Let $G$ be an $M$-regular digraph on the set of vertices $[N]$, $H_1$ be a $D_1$-regular digraph on the set of vertices $[M]$
and $H_2$ be a $D_2$-regular digraph on the set of vertices $[M]$.
Their zig-zag product, denoted by $\z$, is defined to be the graph on the set of vertices $[N]\times[M]$
whose edges are defined as follows: we connect two vertices from $(u,k)$ to $(v,l)$
if it is possible to get from $(u,k)$ to $(v,l)$ by a sequence of the following three steps:
\begin{description}
 \item[Step 1.] Move from $(u,k)$ to $(u,k^\prime)$ where $(k,k^\prime)$ is an edge in $H_1$.
 \item[Step 2.] Move from $(u,k^\prime)$ to $(v,l^\prime)$ where, in $G$, the $k^\prime$-th edge leaving $u$ is the $l^\prime$-th edge entering $v$.
 \item[Step 3.] Move from $(v,l^\prime)$ to $(v,l)$ where $(l^\prime,l)$ is an edge in $H_2$.
\end{description}
Step 1 and 3 change the second component of the current vertex as a transition of the random walk on $H_1$ and $H_2$, respectively.
Step 2 is a deterministic transition using the two-way labelling of $G$.
Therefore their zig-zag product $\z$ is $D_1D_2$-regular.
The two-way labellings of $H_1$ and $H_2$ are used to define the two-way labelling of $\z$ through their rotation maps.
The formal definition of the generalized zig-zag product is as follows.

\begin{Def}[A generalized zig-zag product]
Let $G$ be a two-way labelled $M$-regular graph on the set of vertices
$[N]$ with a rotation map $\Rot_G$,
$H_1$ be a two-way labelled $D_1$-regular graph on the set of vertices
$[M]$ with a rotation map $\Rot_{H_1}$
and $H_2$ be a two-way labelled $D_2$-regular graph on the set of
vertices $[M]$ with a rotation map $\Rot_{H_2}$.
Their zig-zag product $\z$ is defined to be the $D_1D_2$-regular graph
on the set of vertices $[N]\times[M]$ whose rotation map $\Rotz:([N]\times[M])\times([D_1]\times[D_2])
\rightarrow([N]\times[M])\times([D_1]\times[D_2])$ is as follows: $\Rotz((u,k),(i,j))=((v,l),(i^\prime,j^\prime))$
if there exist $k^\prime,l^\prime \in [M]$ such that
\[
(k^\prime,i^\prime) = \Rot_{H_1}(k,i), \quad
(v,l^\prime) = \Rot_G(u,k^\prime),\quad
(l,j^\prime) = \Rot_{H_2}(l^\prime,j).
\]
\end{Def}

Next, we consider a special case of the generalized zig-zag
product where $H_2$ is a trivial 1-regular graph with a single 
self-loop for each vertex.
We denote this trivial graph by $\nullgraph$.  Then we define the
reduced zig-zag product of $G$ and $H$ by
\[
\zp = G\mbox{\textcircled{$\mathrm{z}$}}(H,\nullgraph).
\]
Here $G$ is an $M$-regular digraph on the set of vertices $[N]$ , $H$
is a $D$-regular digraph on the set of vertices $[M]$ and the reduced
zig-zag product $\zp$ is a $D$-regular digraph on the set of vertices
$[N]\times[M]$.  The motivation to consider the reduced zig-zag
product is that the degree of $\zp$ is only $D$ and it still has good
connectivity properties as an expander graph 
as shown below in Theorem \ref{thm:reduced}.

In the reduced zig-zag product, the third step of connecting edges
in the zig-zag product is omitted and the 
edges are defined as follows: we connect two
vertices from $(u,k)$ to $(v,l)$ if it is possible to get from $(u,k)$
to $(v,l)$ by a sequence of the following two steps:
\begin{description}
 \item[Step 1.] Move from $(u,k)$ to $(u,k^\prime)$ where $(k,k^\prime)$ is an edge in $H$.
 \item[Step 2.] Move from $(u,k^\prime)$ to $(v,l)$ where, in $G$, the $k^\prime$-th edge leaving $u$ is the $l^\prime$-th edge entering $v$.
\end{description}
The formal definition is as follows.

\begin{Def}[Reduced zig-zag product]
Let $G$ be a two-way labelled $M$-regular graph on the set of vertices
$[N]$ with a rotation map $\Rot_G$
and $H$ be a two-way labelled $D$-regular graph on the set of vertices
$[M]$ with a rotation map $\Rot_H$.
Their reduced zig-zag product $\zp$ is defined to be the $D$-regular graph
on the set of vertices $[N]\times[M]$ whose rotation map $\Rotzp:([N]\times[M])\times[D]
\rightarrow([N]\times[M])\times[D]$ is as follows: $\Rotzp((u,k),i)=((v,l),j)$
if there exist $k^\prime \in [M]$ such that
\[
(k^\prime,j) = \Rot_{H_1}(k,i),\quad
(v,l) = \Rot_G(u,k^\prime).
\]
\end{Def}

In addition, we define the $t$-th power of a graph, which only replaces the edge
set with the set of all walks of length $t$ in the graph.

\begin{Def}[Powering]
Let $G$ be a two-way labelled $M$-regular graph on the set of vertices
$[N]$ with a rotation map $\Rot_G$.
The $t$-th power of $G$ is the $M^t$-regular graph $G^t$ on the set of vertices $[N]$ whose rotation map
$\Rot_{G^t}:[N]\times[M]^t\rightarrow[N]\times[M]^t$ is defined by
$\Rot_{G^t}(v_0,(k_1,k_2,\dots,k_t)) = (v_t,(l_t,l_{t-1},\dots,l_1))$ where $(v_i,l_i)=\Rot_G(v_{i-1},k_i)$ for $i=1,2,\dots,t$.
\end{Def}

If we denote the transition matrix of $G$ by $A$, then the transition matrix of $G^t$ is $A^t$.

\section{Upper Bounds for Spectral Expansion}
In this section, we derive upper bounds for the spectral expansions of the generalized zig-zag product and the reduced zig-zag product.

\subsection{The Generalized Zig-zag Product}
First, as the  main result of this paper, we derive an upper bound for 
the spectral expansion of the generalized zig-zag product.

\begin{Thm}
\label{thm:main}
If $\lambda(G)\leq\alpha,\, \lambda(H_1)\leq\beta_1$ and $\lambda(H_2)\leq\beta_2$,
then $\lambda(\z)\leq f(\alpha,\beta_1,\beta_2)$, where
\begin{align}
f(\alpha,\beta_1,\beta_2) &=
\frac{1}{2}\left\{\sqrt{\alpha^2(1-\beta_1^2)(1-\beta_2^2)+(\beta_1+\beta_2)^2}
\right. \nonumber\\
& \qquad \  +
\left.
  \sqrt{\alpha^2(1-\beta_1^2)(1-\beta_2^2)+(\beta_1-\beta_2)^2}\right\} .
\label{eq:fabb1}
\end{align}
\end{Thm}

\begin{Rem}
By definition $0\le \alpha, \beta_1, \beta_2 \le 1$ and
\begin{align*}
\alpha^2(1-\beta_1^2)(1-\beta_2^2)+(\beta_1+\beta_2)^2 & =
(1+\beta_1 \beta_2)^2  -(1-\alpha^2)(1-\beta_1^2)(1-\beta_2^2)\\
&\le (1+\beta_1 \beta_2)^2.
\end{align*}
Similarly
\[
\alpha^2(1-\beta_1^2)(1-\beta_2^2)+(\beta_1-\beta_2)^2 
\le (1-\beta_1 \beta_2)^2.
\]
Therefore
\[ f(\alpha,\beta_1,\beta_2) \le \frac{1}{2}
\left\{ \sqrt{(1+\beta_1 \beta_2)^2}+\sqrt{(1-\beta_1 \beta_2)^2} \right\} = 1 \]
with equality holding if
and only if 
\begin{equation}
\label{eq:feq1}
0=(1-\alpha)(1-\beta_1)(1-\beta_2).
\end{equation}
\end{Rem}

The rest of this section is devoted to the proof of (\ref{eq:fabb1}).

\begin{proof}
Let $A,B_1$ and $B_2$ be the transition matrix of the random walk on $G,H_1$ and $H_2$, respectively.
To analyze $\lambda(\z)$, we express $Z$, the transition matrix of $\z$, in terms of $G,H_1$ and $H_2$.
We can decompose $Z$ into the product of three matrices,
corresponding to the three steps in the definition of the edges of $\z$.
Let $\tilde{B}_1$ denote the transition matrix corresponding to the first step.
The first step is only concerned with the
the second component of $[N]\times[M]$.
Hence it is easy to see that 
$\tilde{B}_1=I_N\otimes B_1$, where $I_N$ is the $N\times N$ identity matrix.
Similarly, we have $\tilde{B}_2=I_N\otimes B_2$ where $\tilde{B}_2$ is
the transition matrix corresponding to the third step.
Let $\tilde{A}$ be the transition matrix corresponding to the second step.
Then $\tilde{A}$ is the permutation matrix corresponding to $\Rot_G$, i.e.
$\tilde{A}_{(v,j),(u,i)} = \mathbb{I}\{\Rot_G(u,i)=(v,j)\}$
where $\mathbb{I}\{\cdot\}$ denotes an indicator function which takes
1 if the  condition in the braces is true and 0 otherwise.
Thus $Z$ is written as
\[Z=\tilde{B}_2\tilde{A}\tilde{B}_1 = (I_N\otimes B_2) \, \tilde A \, (I_N
\otimes B_1).
\]

Our aim is to show that $\Vert Z\bm{x}\Vert \leq f(\alpha,\beta_1,\beta_2)\Vert \bm{x}\Vert $ for every $\bm{x}\perp \bm{1}_{MN}$.
In view of the decomposition $Z=\tilde{B}_2\tilde{A}\tilde{B}_1$, we define
\[
\bm{y}=\tilde{B}_1\bm{x}, \ \bm{z}=\tilde{A}\bm{y}, \
\bm{w}=\tilde{B}_2\bm{z}.
\]
For every $u\in [N]$, we define $\bm{x}_u\in\mathbb{R}^M$ by $(x_u)_s = x_{(u,s)}$.
Then, $\bm{x} = \sum_u \bm{e}_u \otimes \bm{x}_u$,
where $\bm{e}_u=(0,\dots,0,1,0,\dots,0)^T$ 
denotes the $u$-th standard basis vector in $\mathbb{R}^N$.
Every $\bm{x}_u$ can be decomposed (uniquely) into $\bm{x}_u = \bm{x}_u^\para + \bm{x}_u^\perp$
where $\bm{x}_u^\para$ is parallel to $\bm{1}_M$ and $\bm{x}_u^\perp$ is orthogonal to $\bm{1}_M$.
Thus, we obtain a decomposition $\bm{x} = \bm{x}^\para + \bm{x}^\perp$
where 
\[
\bm{x}^\para = \sum_u \bm{e}_u \otimes \bm{x}_u^\para, \quad
\bm{x}^\perp = \sum_u \bm{e}_u \otimes \bm{x}_u^\perp.
\]
Since $\bm{x}_u^\para \perp \bm{x}_u^\perp$ for all $u$, we have $\bm{x}^\para \perp \bm{x}^\perp$ and hence
\[
\Vert \bm{x}\Vert ^2 = \Vert \bm{x}^\para\Vert ^2 + \Vert \bm{x}^\perp\Vert ^2.
\]
$\bm{x}^\para$ can also be written as $\bm{x}^\para = \bar{\bm{x}} \otimes \bm{1}_M$,
where $\bar{\bm{x}}\in\mathbb{R}^N$ is defined by
$\bar{x}_u = (1/M)\sum_{s\in[M]}x_{(u,s)}$.
Since $\bm{x}$ and $\bm{x}^\perp$ are both orthogonal to $\bm{1}_{MN}$, so is $\bm{x}^\para$
and hence also $\bar{\bm{x}}$ is orthogonal to $\bm{1}_N$.
We decompose $\bm{y},\bm{z}$ and $\bm{w}$ in the same way as $\bm{x}$.

We now show several relations among $\bm{x},\bm{y},\bm{z}$ and $\bm{w}$.
First we consider the relation between  $\bm{x}$ and $\bm{y}$. Since
\begin{eqnarray*}
\bm{y} &=& \tilde{B}_1\bm{x} \\
&=& \tilde{B}_1\bm{x}^\para + \tilde{B}_1\bm{x}^\perp \\
&=& \sum_{u\in[N]}\bm{e}_u\otimes B_1\bm{x}_u^\para + \sum_{u\in[N]}\bm{e}_u\otimes B_1\bm{x}_u^\perp 
\end{eqnarray*}
and $B_1\bm{x}_u^\para=\bm{x}_u^\para\parallel\bm{1}_M,B_1\bm{x}_u^\perp\perp\bm{1}_M$ for every $u$,
we have $\bm{y}^\para=\bm{x}^\para$ and $\bm{y}^\perp=\tilde{B}_1\bm{x}^\perp$.
Note that $\Vert \bm{y}_u^\perp\Vert =\Vert B_1\bm{x}_u^\perp\Vert \leq\beta_1\Vert \bm{x}_u^\perp\Vert $ for every $u\in[M]$.
Thus, we have
\begin{equation}
\label{eq:yx}
\Vert \bm{y}^\para\Vert =\Vert \bm{x}^\para\Vert ,\quad 
\Vert \bm{y}^\perp\Vert \leq\beta_1\Vert \bm{x}^\perp\Vert .
\end{equation}
The relation between $\bm{z}$ and $\bm{w}$ is similar and we have
\begin{equation}
\label{eq:wz}
\Vert \bm{w}^\para\Vert =\Vert \bm{z}^\para\Vert ,\quad
\Vert \bm{w}^\perp\Vert \leq\beta_2\Vert \bm{z}^\perp\Vert .
\end{equation}

Now consider the relation between 
$\bm{y}$ and $\bm{z}$.  Since $\tilde{A}$ is a permutation matrix, we have
\begin{equation}
\label{eq:zy}
\Vert \bm{z}\Vert =\Vert \tilde{A}\bm{y}\Vert =\Vert \bm{y}\Vert .
\end{equation}
Furthermore for $v\in[N]$
\begin{eqnarray*}
 (\overline{\tilde{A}\bm{y}^\para})_v &=& \frac{1}{M}\sum_{j\in[M]}(\tilde{A}\bm{y}^\para)_{(v,j)} \\
 &=& \frac{1}{M}\sum_{j\in[M]}\sum_{(u,i)\in[N]\times[M]} \tilde{A}_{(v,j),(u,i)}y^\para_{(u,i)} \\
 &=& \frac{1}{M}\sum_{u\in[N]}\sum_{i,j\in[M]} \mathbb{I}\{\mathrm{Rot}_G(u,i)=(v,j)\}\bar{y}_u \\
 &=& \sum_{u\in[N]}\frac{|\mbox{the number of edges from $u$ to $v$}|}{M}\cdot \bar{y}_u \\
 &=& \sum_{u\in[N]}A_{vu}\bar{y}_u \\
 &=& (A\bar{\bm{y}})_v,
\end{eqnarray*}
and $\Vert A\bar{\bm{y}}\Vert \leq\alpha\Vert \bar{\bm{y}}\Vert $.  Therefore we have
\begin{equation}
\label{eq:zypara}
\Vert (\tilde{A}\bm{y}^\para)^\para\Vert =\Vert (\overline{\tilde{A}\bm{y}^\para})\otimes\bm{1}_M\Vert 
=\Vert A\bar{\bm{y}}\otimes\bm{1}_M\Vert \leq\alpha\Vert
\bar{\bm{y}}\otimes\bm{1}_M\Vert =\alpha\Vert \bm{y}^\para\Vert  .
\end{equation}

We will now prove 
$\Vert \bm{w}\Vert \leq f(\alpha,\beta_1,\beta_2)\Vert \bm{x}\Vert $ considering two cases,
depending on the size of the norm $\Vert \bm{x}^\para\Vert $.

\smallskip
\noindent
\textbf{Case 1:} $\Vert \bm{x}^\para\Vert \geq\alpha\Vert \bm{x}\Vert $

Since $\Vert \bm{y}^\para\Vert  = \Vert \bm{x}^\para\Vert  \geq \alpha\Vert \bm{x}\Vert  \geq \alpha\Vert \bm{y}\Vert $
and hence $\Vert \bm{y}^\perp\Vert \leq\sqrt{1-\alpha^2}\Vert \bm{y}\Vert $, we have
\[ 1 - \frac{\alpha}{\sqrt{1-\alpha^2}}\frac{\Vert \bm{y}^\perp\Vert }{\Vert \bm{y}^\para\Vert } \geq 0. \]
Using this, the triangle inequality yields
\begin{eqnarray}
\Vert \bm{z}^\para\Vert  &=& \Vert (\tilde{A}\bm{y})^\para\Vert  \nonumber\\
&=& \Vert \{\tilde{A}(\bm{y}^\para+\bm{y}^\perp)\}^\para\Vert  \nonumber\\
&=& \left\Vert  \left\{\tilde{A}\cdot\left(1-\frac{\alpha}{\sqrt{1-\alpha^2}}
\frac{\Vert \bm{y}^\perp\Vert }{\Vert \bm{y}^\para\Vert }\right)
\bm{y}^\para\right\}^\para+\left\{\tilde{A}\left(\frac{\alpha}{\sqrt{1-\alpha^2}}
\frac{\Vert \bm{y}^\perp\Vert }{\Vert \bm{y}^\para\Vert }\bm{y}^\para+\bm{y}^\perp\right)\right\}^\para\right\Vert \nonumber\\
&\leq& \left(1-\frac{\alpha}{\sqrt{1-\alpha^2}}\frac{\Vert \bm{y}^\perp\Vert }{\Vert \bm{y}^\para\Vert }\right)
\Vert (\tilde{A}\bm{y}^\para)^\para\Vert +\left\Vert \left\{\tilde{A}\left(\frac{\alpha}{\sqrt{1-\alpha^2}}
\frac{\Vert \bm{y}^\perp\Vert }{\Vert \bm{y}^\para\Vert }\bm{y}^\para+\bm{y}^\perp\right)\right\}^\para\right\Vert .
\label{eq:case1zy1}
\end{eqnarray}
From (\ref{eq:zypara}) the first term of the right hand of (\ref{eq:case1zy1}) is bounded by
\[ \left(1-\frac{\alpha}{\sqrt{1-\alpha^2}}\frac{\Vert \bm{y}^\perp\Vert }{\Vert \bm{y}^\para\Vert }\right)
\Vert (\tilde{A}\bm{y}^\para)^\para\Vert  \leq
\left(1-\frac{\alpha}{\sqrt{1-\alpha^2}}\frac{\Vert \bm{y}^\perp\Vert }{\Vert \bm{y}^\para\Vert }\right)
\cdot \alpha\Vert \bm{y}\Vert  \]
and from (\ref{eq:zy}) the second term is bounded by
\begin{eqnarray*}
\left\Vert \left\{\tilde{A}\left(\frac{\alpha}{\sqrt{1-\alpha^2}}
\frac{\Vert \bm{y}^\perp\Vert }{\Vert \bm{y}^\para\Vert }\bm{y}^\para+\bm{y}^\perp\right)\right\}^\para\right\Vert 
&\leq& \left\Vert \tilde{A}\left(\frac{\alpha}{\sqrt{1-\alpha^2}}
\frac{\Vert \bm{y}^\perp\Vert }{\Vert \bm{y}^\para\Vert }\bm{y}^\para+\bm{y}^\perp\right)\right\Vert  \\
&=& \left\Vert \frac{\alpha}{\sqrt{1-\alpha^2}}
\frac{\Vert \bm{y}^\perp\Vert }{\Vert \bm{y}^\para\Vert }\bm{y}^\para+\bm{y}^\perp\right\Vert .
\end{eqnarray*}
Therefore we obtain
\begin{eqnarray}
\Vert \bm{z}^\para\Vert 
&\leq& \left(1-\frac{\alpha}{\sqrt{1-\alpha^2}}\frac{\Vert \bm{y}^\perp\Vert }{\Vert \bm{y}^\para\Vert }\right)\cdot
\alpha\Vert \bm{y}^\para\Vert +\left\Vert \frac{\alpha}{\sqrt{1-\alpha^2}}\frac{\Vert \bm{y}^\perp\Vert }{\Vert \bm{y}^\para\Vert }
\bm{y}^\para+\bm{y}^\perp\right\Vert  \nonumber\\
&=& \alpha\Vert \bm{y}^\para\Vert -\frac{\alpha^2}{\sqrt{1-\alpha^2}}\Vert \bm{y}^\perp\Vert +
\sqrt{\left(\frac{\alpha}{\sqrt{1-\alpha^2}}\frac{\Vert \bm{y}^\perp\Vert }{\Vert \bm{y}^\para\Vert }\Vert \bm{y}^\para\Vert \right)^2
+\Vert \bm{y}^\perp\Vert ^2} \nonumber\\
&=& \alpha\Vert \bm{y}^\para\Vert -\frac{\alpha^2}{\sqrt{1-\alpha^2}}\Vert \bm{y}^\perp\Vert +
\frac{1}{\sqrt{1-\alpha^2}}\Vert \bm{y}^\perp\Vert  \nonumber\\
&=& \alpha\Vert \bm{y}^\para\Vert  + \sqrt{1-\alpha^2}\Vert \bm{y}^\perp\Vert .
\label{eq:case1zy2}
\end{eqnarray}
Thus, from (\ref{eq:yx}), (\ref{eq:wz}), (\ref{eq:zy}) and (\ref{eq:case1zy2}), $\Vert \bm{w}\Vert $ is bounded as follows:
\begin{eqnarray}
\Vert \bm{w}\Vert ^2 &=& \Vert \bm{w}^\para\Vert ^2 + \Vert \bm{w}^\perp\Vert ^2 \nonumber\\
&\leq& \Vert \bm{z}^\para\Vert ^2 + \beta_2^2\Vert \bm{z}^\perp\Vert ^2 \nonumber\\
&=& (1-\beta_2^2)\Vert \bm{z}^\para\Vert ^2 + \beta_2^2\Vert \bm{z}\Vert ^2 \nonumber\\
&\leq& (1-\beta_2^2)(\alpha\Vert \bm{y}^\para\Vert +\sqrt{1-\alpha^2}\Vert \bm{y}^\perp\Vert )^2 + \beta_2^2\Vert \bm{y}\Vert ^2 \nonumber\\
&\leq& (1-\beta_2^2)(\alpha\Vert \bm{x}^\para\Vert +\sqrt{1-\alpha^2}\cdot\beta_1\Vert \bm{x}^\perp\Vert )^2
+ \beta_2^2(\Vert \bm{x}^\para\Vert ^2+\beta_1^2\Vert \bm{x}^\perp\Vert ^2) \nonumber\\
&=&(\alpha^2+\beta_2^2-\alpha^2\beta_2^2)\Vert \bm{x}^\para\Vert ^2+
(\beta_1^2-\alpha^2\beta_1^2+\alpha^2\beta_1^2\beta_2^2)\Vert \bm{x}^\perp\Vert ^2 \nonumber\\
&& + 2\alpha\beta_1\sqrt{1-\alpha^2}(1-\beta_2^2)\Vert \bm{x}^\para\Vert \cdot\Vert \bm{x}^\perp\Vert .
\label{eq:wx}
\end{eqnarray}
Now, it is straightforward to maximize the right hand side of (\ref{eq:wx}) 
subject to $\Vert \bm{x}^\para\Vert ^2 + \Vert \bm{x}^\perp\Vert ^2 =
\Vert \bm{x}\Vert ^2$ by the Lagrange multiplier method 
and we obtain
\begin{align*}
\Vert \bm{w}\Vert ^2 &\leq \left. \frac{1}{2} \right[ \alpha^2(1-\beta_1^2)(1-\beta_2^2) + (\beta_1^2+\beta_2^2)  \\
& \left. + \sqrt{\{\alpha^2(1-\beta_1^2)(1-\beta_2^2)+(\beta_1+\beta_2)^2\}\cdot
\{\alpha^2(1-\beta_1^2)(1-\beta_2^2)+(\beta_1-\beta_2)^2\} }\right] \cdot\Vert \bm{x}\Vert ^2 \\
&=
\frac{1}{4}\left\{\sqrt{\alpha^2(1-\beta_1^2)(1-\beta_2^2)+(\beta_1+\beta_2)^2}
 \right.  \\
& \qquad\quad  \left. +\sqrt{\alpha^2(1-\beta_1^2)(1-\beta_2^2)+(\beta_1-\beta_2)^2}\right\}^2\cdot\Vert \bm{x}\Vert ^2 \\
&= f(\alpha,\beta_1,\beta_2)^2\Vert \bm{x}\Vert ^2.
\end{align*}

\noindent\textbf{Case 2:} $\Vert \bm{x}^\para\Vert <\alpha\Vert \bm{x}\Vert $

From (\ref{eq:wz}) and (\ref{eq:zy}), we have
\begin{eqnarray*}
\Vert \bm{w}\Vert ^2 &\leq& \Vert \bm{z}\Vert ^2 = \Vert \bm{y}\Vert ^2 \\
&=& \Vert \bm{y}^\para\Vert ^2 + \Vert \bm{y}^\perp\Vert ^2 \\
&\leq& \Vert \bm{x}^\para\Vert ^2 + \beta_1^2\Vert \bm{x}^\perp\Vert ^2 \\
&=& (1-\beta_1^2)\Vert \bm{x}^\para\Vert ^2 + \beta_1^2\Vert \bm{x}\Vert ^2 \\
&<& \{(1-\beta_1^2)\alpha^2 + \beta_1^2\}\Vert \bm{x}\Vert ^2
\end{eqnarray*}
and this is smaller than $f(\alpha,\beta_1,\beta_2)^2\Vert \bm{x}\Vert ^2$.

From the above two cases we conclude that $\Vert \bm{w}\Vert =\Vert Z\bm{x}\Vert \leq
f(\alpha,\beta_1,\beta_2)\Vert \bm{x}\Vert $ and 
\[ \lambda(\z) = \max_{\bm{x}\perp\bm{1}_{MN}}\frac{\Vert Z\bm{x}\Vert }{\Vert \bm{x}\Vert } \leq f(\alpha,\beta_1,\beta_2). \]

\end{proof}

\subsection{The Reduced Zig-zag Product}
Here we derive an upper bound for the spectral expansion of the powers
of the reduced zig-zag product.  Note that (\ref{eq:fabb1}) in Theorem
\ref{thm:main} is not useful for the reduced zig-zag product because
$\lambda(\nullgraph)=1$ and $f(\alpha,\beta_1,1)=1$.
In fact $\lambda(\zp) = 1$.  However the spectral
expansion of higher  powers of  the reduced zig-zag product behaves as
in the generalized zig-zag product as shown in the following theorem.

\begin{Thm}
\label{thm:reduced}
$\lambda(\zp) = 1$.
If $\lambda(G)\leq\alpha$ and $\lambda(H)\leq\beta$, then
\begin{equation}
\label{eq:reduced-power}
\lambda((\zp)^k) \leq f^\prime(\alpha,\beta)^{k-1}, \quad k=2,3,\ldots
\end{equation}
where
\[ f^\prime(\alpha,\beta)=f(\alpha,\sqrt{\beta},\sqrt{\beta})=
\frac{1}{2}\alpha(1-\beta)+\frac{1}{2}\sqrt{\alpha^2(1-\beta)^2+4\beta} . \]
\end{Thm}

\begin{proof}
Let $Z^\prime$ be the transition matrix of the random walk on $\zp$.
In the same way as the proof of Theorem \ref{thm:main}, $Z^\prime$ is decomposed into $Z^\prime = \tilde{A}\tilde{B}$.
When $x=x^\para$, $\tilde{B}x=x$ and hence $\Vert Z^\prime x\Vert  = \Vert \tilde{A}\tilde{B}x\Vert  = \Vert \tilde{A}x\Vert  = \Vert x\Vert $. 
This implies that  $\lambda(\zp)=1$.

To prove 
(\ref{eq:reduced-power}),  we consider a singular value decomposition: $B = P\Sigma Q^T$
where $P$ and $Q$ are orthogonal matrices and $\Sigma$ is a diagonal
matrix of the singular values of $B$.
We denote $P=(\bm{p}_1,\dots,\bm{p}_M),Q=(\bm{q}_1,\dots,\bm{q}_M),\Sigma=diag(\sigma_1,\dots,\sigma_M)$
and set $\sigma_1=1,\bm{p}_1=\bm{q}_1=\bm{1}_M$.
Then, since $P^TP=I_M$, we have
\begin{eqnarray*}
B &=& P\Sigma Q^T \\
&=& P\sqrt{\Sigma}P^T P\sqrt{\Sigma}Q^T \\
&=:& B_1 B_2
\end{eqnarray*}
where $B_1=P\sqrt{\Sigma}P^T, B_2=P\sqrt{\Sigma}Q^T$.
Then, $\tilde{B} = \tilde{B}_1\tilde{B}_2$ where $\tilde{B}_1=I_N \otimes B_1, \tilde{B}_2=I_N \otimes B_2$,
and for $k=2,3,\dots$, we can decompose ${Z^\prime}^k$, the transition matrix of $(\zp)^k$,
into
\begin{eqnarray*}
{Z^\prime}^k &=& (\tilde{A}\tilde{B})^k \\
&=& (\tilde{A}\tilde{B}_1\tilde{B}_2)^k \\
&=& \tilde{A}\tilde{B}_1(\tilde{B}_2\tilde{A}\tilde{B}_1)^{k-1}\tilde{B}_2.
\end{eqnarray*}
Since the singular values of both $B_1$ and $B_2$ are the square roots of those of $B$,
the second largest singular values of $B_1$ and $B_2$ are both smaller
than or equal to $\sqrt{\beta}$.
Also, since the left-singular and right-singular vectors of $B_1,B_2$
corresponding to the singular value $\sqrt{\sigma_1}=1$ are all $\bm{p}_1=\bm{q}_1=\bm{1}_M$,
we have $B_1\bm{1}_M=B_2\bm{1}_M=\bm{1}_M$.
It follows that $\tilde{B}_1$ and $\tilde{B}_2$ have the same
properties as in the proof of Theorem \ref{thm:main} when
$\sqrt{\beta}$ is substituted for $\beta_1$ and $\beta_2$.
Therefore we have $\Vert \tilde{B}_2\tilde{A}\tilde{B}_1\bm{x}\Vert \leq f(\alpha,\sqrt{\beta},\sqrt{\beta})\Vert \bm{x}\Vert 
= f^\prime(\alpha,\beta)\Vert \bm{x}\Vert $ for every $\bm{x}\perp \bm{1}_{MN}$.
Now 
\begin{eqnarray*}
\Vert {Z^\prime}^k\bm{x}\Vert  &=& \Vert \tilde{A}\tilde{B}_1(\tilde{B}_2\tilde{A}\tilde{B}_1)^{k-1}(\tilde{B}_2\bm{x})\Vert  \\
&\leq& \Vert (\tilde{B}_2\tilde{A}\tilde{B}_1)^{k-1}(\tilde{B}_2\bm{x})\Vert  \\
&\leq& f^\prime(\alpha,\beta)^{k-1}\Vert \tilde{B}_2\bm{x}\Vert  \\
&\leq& f^\prime(\alpha,\beta)^{k-1}\Vert \bm{x}\Vert 
\end{eqnarray*}
and hence 
\[ \lambda((\zp)^k)  =  \max_{\bm{x}\perp\bm{1}_{MN}}
\frac{\Vert {Z^\prime}^k\bm{x}\Vert }{\Vert \bm{x}\Vert } \leq f^\prime(\alpha,\beta)^{k-1}. \]

\end{proof}

\section{Numerical experiments}
In this section, we compare spectral expansions for zig-zag products
with their bounds given in Theorem \ref{thm:main}
and Theorem \ref{thm:reduced}.
We use a random $M$-regular digraph $G$ on the set of vertices $[N]$ and $D$-regular digraph $H$ on the set of vertices $[M]$
in three cases:\\[8pt]
(i) $N=50,M=40,D=30.$ \\[4pt]
(ii) $N=30,M=20,D=10.$ \\[4pt]
(iii) $N=10,M=5,D=3.$   \\[8pt]
The random graphs are generated by the following algorithm based on the configuration model (\cite{Novi}, \cite{Wor}).

\begin{Alg}[configulation model]
In generating a random $M$-regular digraph $G=(V,E)$ on the set of
vertices $[N]$, we take the following steps
(generating $H$ is similar):
\begin{description}
 \item[Step 1.] Create two vectors $\bm{x}$ and $\bm{y}$, each a random permutation of the integers from 1 to $MN$.
 \item[Step 2.] Reassign all the entries in the vectors with their values mod $N$ (integers from 1 to $N$).
 \item[Step 3.] Construct the graph $G$ by defining vertex set $V(G)=[N]$ and edge set $E(G)=\{(x_i,y_i):i=1,\dots,MN\}$. 
\end{description}
\end{Alg}

Also we generate 100 random  two-way labellings for each $G$ in order to define the zig-zag products of $G$ and $H$.
We do not need a two-way labelling for $H$ here, since it is used to define
only the two-way labelling for the zig-zag products, which is irrelevant to the spectral expansion of the products.
Thus, we obtain 100 zig-zag products $G\mbox{\textcircled{$\mathrm{z}$}}H=G\mbox{\textcircled{$\mathrm{z}$}}(H,H)$
and reduced zig-zag products $\zp$ in each case.

First, we compare the spectral expansion for the zig-zag product
$\lambda=\lambda(G\mbox{\textcircled{$\mathrm{z}$}}H)$ with their bound
$f=f(\lambda(G),\lambda(H),\lambda(H))$ given in Theorem \ref{thm:main}.
The spectral expansion is computed as the second largest singular value of the transition matrix of the graph.
In each case, we computed the average and the maximum of $\lambda$.
The results are shown in Table \ref{tab:zigzag}.
The gap between $f$ and the maximum of $\lambda$ is about 0.1 in each case
and hence the bound $f$ is tight to some degree.
Nevertheless, the gap is much larger than the range of $\lambda$ in the case (i) and (ii).
We infer that the singular values of the transition matrices of $G$ and $H$ smaller than $\lambda(G)$ and $\lambda(H)$,
respectively, reduce $\lambda(G\mbox{\textcircled{$\mathrm{z}$}}H)$ from our bound $f(\lambda(G),\lambda(H),\lambda(H))$.

\begin{table}[h] 
\begin{center}
\begin{tabular}{|c|c|c|c|c|c|c|} \hline
& $\lambda(G)$ & $\lambda(H)$ & ave $\lambda$ & $\max \lambda$ & $f$ \\ \hline
(i)&0.2931153 &0.3334984 &0.3692421 &0.3708440 &0.4882911 \\ \hline
(ii)&0.4184724 &0.5226591 &0.5522197 &0.5606170 &0.6964135 \\ \hline
(iii)&0.5909580 &0.8047379 &0.8294209 &0.8610790 &0.9155723 \\ \hline
\end{tabular}
\caption{Comparison of $\lambda=\lambda(G\mbox{\textcircled{$\mathrm{z}$}}H)$
  and its upper bound $f=f(\lambda(G),\lambda(H),\lambda(H))$}
\label{tab:zigzag}
\end{center}
\end{table}

\begin{figure}[p] 
\begin{center}
\includegraphics*[scale=0.8,angle=-90]{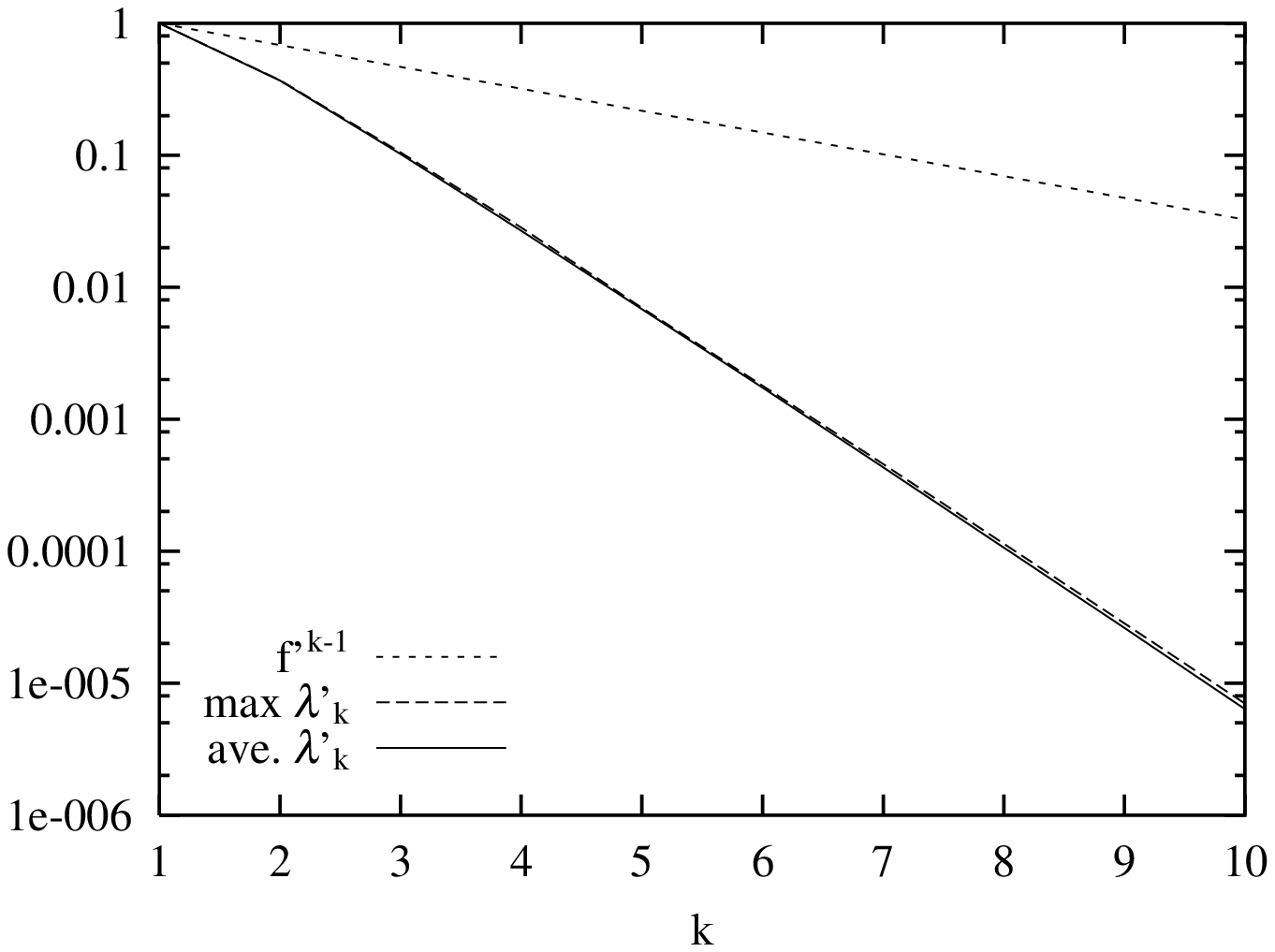}
\end{center}
\caption{Comparison of
  $\lambda^\prime_k=\lambda((\zp)^k)$
  and its upper bound ${f^\prime}^{k-1}$ $(N=50,M=40,D=30)$}
\label{fig:reduced1}
\end{figure}

\begin{figure}[p] 
\begin{center}
\includegraphics*[scale=0.8,angle=-90]{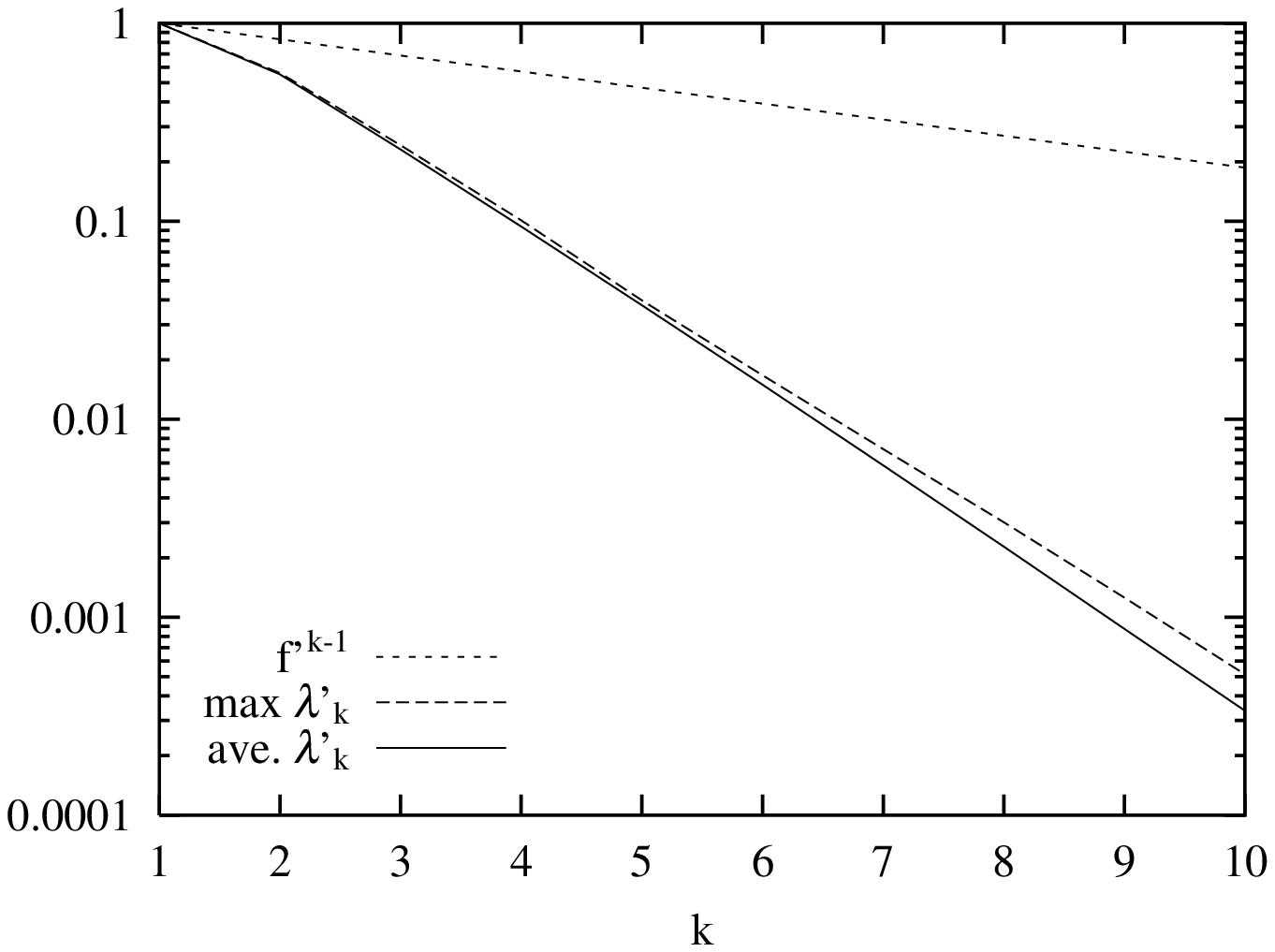}
\end{center}
\caption{Comparison of
  $\lambda^\prime_k=\lambda((\zp)^k)$
  and its upper bound ${f^\prime}^{k-1}$ $(N=30,M=20,D=10)$}
\label{fig:reduced2}
\end{figure}

\begin{figure}[p] 
\begin{center}
\includegraphics*[scale=0.8,angle=-90]{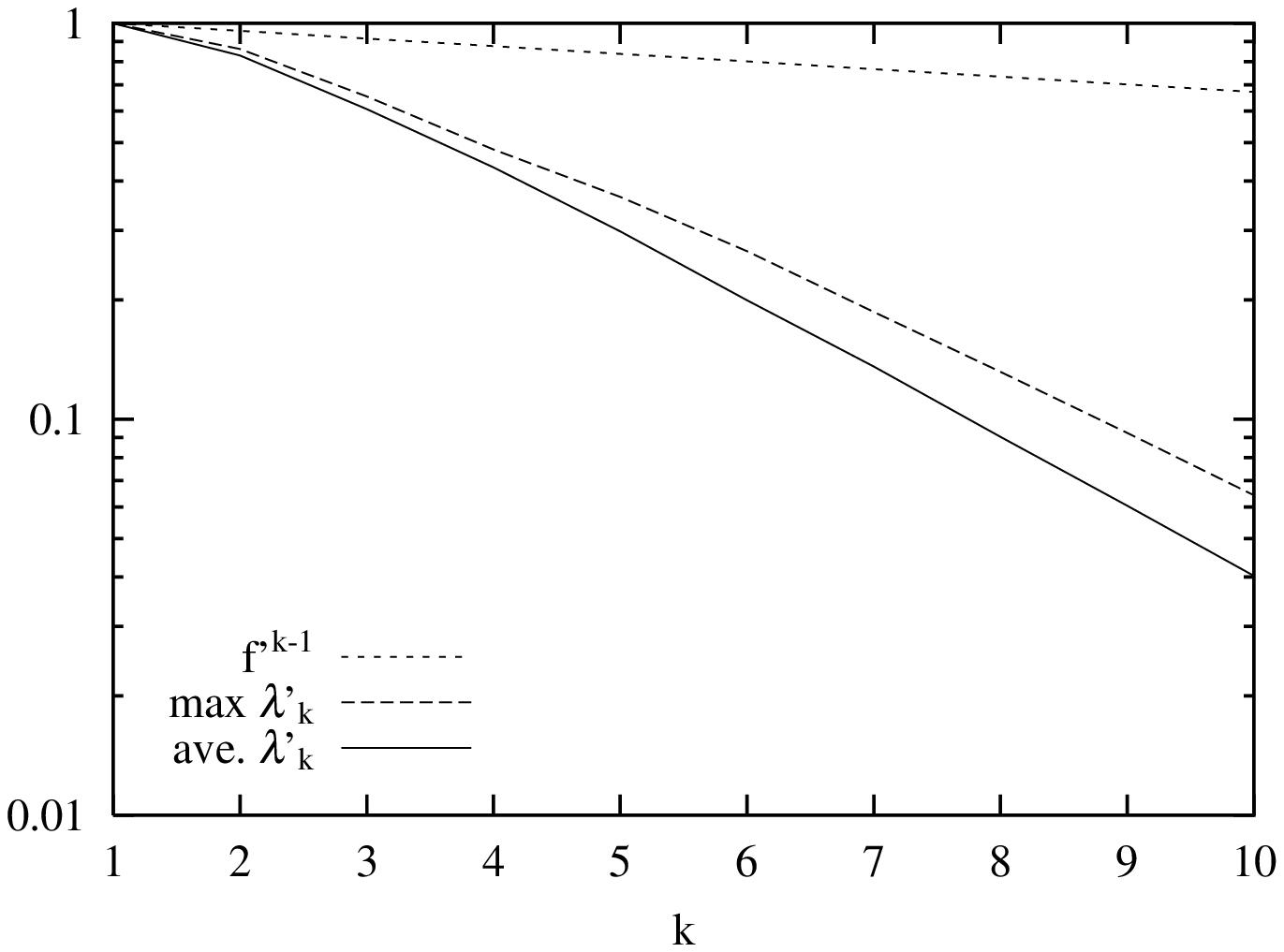}
\end{center}
\caption{Comparison of
  $\lambda^\prime_k=\lambda((\zp)^k)$
  and its upper bound ${f^\prime}^{k-1}$ $(N=10,M=5,D=3)$}
\label{fig:reduced3}
\end{figure}

Next, for $k=1,\dots,10$ we compare the spectral expansion for the $k$-th power of the reduced zig-zag product
$\lambda_k^\prime=\lambda((\zp)^k)$ with their bounds ${f^\prime}^{k-1}=f^\prime(\lambda(G),\lambda(H))^{k-1}$
given in Theorem \ref{thm:reduced}.
We show the results in Table \ref{tab:reduced1}-\ref{tab:reduced3}
and the graphs of the results in Figure \ref{fig:reduced1}-\ref{fig:reduced3}.
The vertical axis is the logarithm of each variables and the horizontal axis is the degree of power $k$.
We can see that the rate at which $\lambda^\prime_k$ decreases as $k$ increases
is much smaller than $f^\prime$.
It is because $\lambda^\prime_k$ decreases asymptotically at the rate of the second largest eigenvalue (in absolute value) of $Z^\prime$,
which is smaller than the second largest singular value of $\tilde{B}_2\tilde{A}\tilde{B}_1$
bounded by $f^\prime$ in the proof of Theorem \ref{thm:reduced}.

\begin{table}[p] 
\begin{center}
\begin{tabular}{|c|c|c|c|} \hline
$k$ &ave $\lambda^\prime_k$ &$\max \lambda^\prime_k$ &${f^\prime}^{k-1}$ \\ \hline
1 &1 &1 &1 \\ \hline
2 &0.3692421 &0.3708440 &0.6833770\\ \hline
3 &0.1022944 &0.1053438 &0.4670042\\ \hline
4 &0.0268145 &0.0284278 &0.3191399\\ \hline
5 &0.0068423 &0.0069972 &0.2180929\\ \hline
6 &0.0017359 &0.0017915 &0.1490397\\ \hline
7 &0.0004315 &0.0004564 &0.1018503\\ \hline
8 &0.0001064 &0.0001137 &0.0696022\\ \hline
9 &0.0000263 &0.0000284 &0.0475645\\ \hline
10 &0.0000064 &0.0000070 &0.0325045\\ \hline
\end{tabular}
\caption{Comparison of
  $\lambda^\prime_k=\lambda((\zp)^k)$
  and its upper bound ${f^\prime}^{k-1}$ $(N=50,M=40,D=30)$}
\label{tab:reduced1}
\end{center}
\end{table}

\begin{table}[p] 
\begin{center}
\begin{tabular}{|c|c|c|c|} \hline
$k$ & ave $\lambda^\prime_k$ &$\max \lambda^\prime_k$ &${f^\prime}^{k-1}$ \\ \hline
1& 1& 1& 1\\ \hline
2& 0.5522197& 0.5606170& 0.8296951\\ \hline
3& 0.2304273& 0.2419442& 0.6883940\\ \hline
4& 0.0942567& 0.1013095& 0.5711572\\ \hline
5& 0.0377101& 0.0399611& 0.4738863\\ \hline
6& 0.0149260& 0.0166791& 0.3931812\\ \hline
7& 0.0058516& 0.0070592& 0.3262205\\ \hline
8& 0.0022697& 0.0030180& 0.2706636\\ \hline
9& 0.0008743& 0.0012554& 0.2245682\\ \hline
10& 0.0003360& 0.0005136& 0.1863232\\ \hline
\end{tabular}
\caption{Comparison of
  $\lambda^\prime_k=\lambda((\zp)^k)$
  and its upper bound ${f^\prime}^{k-1}$  $(N=30,M=20,D=10)$}
\label{tab:reduced2}
\end{center}
\end{table}

\begin{table}[p] 
\begin{center}
\begin{tabular}{|c|c|c|c|} \hline
$k$ & ave $\lambda^\prime_k$ &$\max \lambda^\prime_k$ &${f^\prime}^{k-1}$ \\ \hline
1 &1 &1 &1\\ \hline
2 &0.8294209 &0.8610790 &0.9566212\\ \hline
3 &0.6072055 &0.6531134 &0.9151240\\ \hline
4 &0.4327482 &0.4799104 &0.8754270\\ \hline
5 &0.2983260 &0.3643133 &0.8374520\\ \hline
6 &0.1998155 &0.2659929 &0.8011243\\ \hline
7 &0.1357641 &0.1866002 &0.7663725\\ \hline
8 &0.0902645 &0.1317873 &0.7331281\\ \hline
9 &0.0604630 &0.0923539 &0.7013259\\ \hline
10 &0.0402301 &0.0642308 &0.6709032\\ \hline
\end{tabular}
\caption{Comparison of
  $\lambda^\prime_k=\lambda((\zp)^k)$
  and its upper bound ${f^\prime}^{k-1}$ $(N=10,M=5,D=3)$}
\label{tab:reduced3}
\end{center}
\end{table}

\end{document}